\documentclass[a4paper,12pt]{article}
\usepackage{amsmath,amssymb} \usepackage{euler,eucal}
\usepackage[utf8]{inputenc}
\usepackage[T1]{fontenc}
\usepackage[english]{babel}
\usepackage{tgpagella}
\usepackage{braket}
\usepackage{accents}
\usepackage[unicode=true,plainpages=false,breaklinks=true]{hyperref}
\usepackage{graphicx}
\hypersetup{colorlinks=false,pdfborder={0 0 0.1},linkcolor=magenta,anchorcolor=magenta,urlcolor=blue,citecolor=blue}

\def\email#1{\texttt{#1}}
\usepackage[left=30mm,right=30mm,top=20mm,bottom=30mm,a4paper]{geometry}
\usepackage{algorithmic}
\usepackage[ruled]{algorithm}
\usepackage{amsthm}
\theoremstyle{definition}

\usepackage{tikz} 
\usepackage{pgfplots} \usetikzlibrary{plotmarks}
\let\le=\leqslant

\let\leq=\leqslant

\def\O{\mathcal{O}}  

\def\hat{\mathaccent"705E }
\def\aprod#1#2{\Braket{#1 | #2}}
\def\Id{I}
\def\H{\hat H}
\def\hH{\hat{\hat H}}
\def\Q{Q}
\def\hR{\hat{\hat R}}
\def\s{\hat\sigma}

\def\lam{\lambda}
\def\d{\partial}

\def\ds{\mathrm{d}\,s}
\def\phi{\varphi}
\def\i{i}
\def\e{\mathrm{1}}

\allowdisplaybreaks
\begin{document}

\title{Tensor product approach to quantum control\thanks{
 This work is supported by EPSRC grants EP/P033954/1 (D.Q. and D.S.) and EP/M019004/1 (S.D.).
}}
\author{Diego Qui\~nones Valles\footnotemark[2], Sergey Dolgov\footnotemark[3] and Dmitry Savostyanov\footnotemark[4]}
\date{28 February 2019}
\maketitle

\renewcommand{\thefootnote}{\fnsymbol{footnote}}
\footnotetext[2]{University of Brighton, UK,    \email{d.quinonesvalles@brighton.ac.uk}}
\footnotetext[3]{University of Bath, UK,        \email{s.dolgov@bath.ac.uk}}
\footnotetext[4]{University of Brighton, UK,    \email{d.savostyanov@brighton.ac.uk}}

\begin{abstract}
 In this proof-of-concept paper we show that tensor product approach is efficient for control of large quantum  systems, such as Heisenberg spin wires, which are essential for emerging quantum computing technologies.
 We compute optimal control sequences using GRAPE method, 
  applying the recently developed tAMEn algorithm to calculate evolution of quantum states represented in the tensor train format to reduce storage.
 Using tensor product algorithms we can overcome the curse of dimensionality and compute the optimal control pulse for a $41$ spin system on a single workstation with fully controlled accuracy and huge savings of computational time and memory.
 The use of tensor product algorithms opens new approaches for development of quantum computers with $50$ to $100$ qubits.
\end{abstract}

\section{Introduction}\label{sec:intro}                     
 Quantum control plays a central part in both established technologies such as nuclear magnetic resonance (NMR) and emerging technologies such as quantum computing~\cite{quaint-2015}.
 The `hardware' of a quantum computer is a quantum system (e.g. a spin system), which we manipulate with time--dependent fields (e.g. laser pulses).
 Our goal is to design a pulse sequence that steers the system from an initial state $\ket{\psi_0}$ to a desired state $\ket{\psi_T}$ during time~$T$ --- a `firmware', so to say.
 The main challenge for numerical methods is the notoriously known \emph{curse of dimensionality} --- the dimension of the state space growing exponentially with the number of particles $d.$
 Even for simplest spin--$\frac12$ systems, the wavefunction $\ket{\psi},$ the density matrix $\rho,$ and the Hamiltonian $\H$ require at least $\O(2^d)$ storage.
 To find the optimal trajectory we need to compute $\rho(t)$ in $100$s to $1000$s points of time.
 This leads to gigabyte--scale storage for $d \simeq 20$ spins and becomes unfeasible for $d \gtrsim 30.$

 The main approaches used currently are
  the sparse format for $\H$, $\ket{\psi}$'s, $\rho$'s (only large entries are stored, zeroes and small elements are discarded), 
  and dimensionality reduction --- the equations of system dynamics are projected into a subspace of low--lying eigenstates~\cite{spinwire-simple-2005} or its heuristically chosen equivalent~\cite{kuprov-2007,crab-2011}.

 Tensor product algorithms are a novel approach to high--dimen\-si\-onal problems.
 Based on the idea of separation of variables, they approximate high--dimen\-si\-on\-al data with a tensor product of low--dimensional factors.
 For example, when interaction of particles is a superposition of a few near--neighbour interactions, the corresponding quantum state is a superposition of a few unentangled states, i.e it has a low Schmidt rank. 
 The corresponding tensor format also has a low rank, providing good compression and low storage costs.
 By maintaining the compressed format for all operations we can compute the system's trajectory and optimise it without any uncontrollable truncations or heuristic assumptions.
 Model examples of such systems include Ising and Heisenberg spin chains~\cite{schollwock-2011},
 more realistic applications include quantum wires~\cite{jeckelmann-wire-2017} and backbones of simple protein molecules~\cite{sdwk-nmr-2014}.

 In quantum physics a tensor product format was 
   first proposed in 1970s as the renormalization group formalism~\cite{wilson-nrg-1975}, 
   followed in 1990s by MPS/MPO~\cite{fannes-mps-1992,klumper-mps-1993} and DMRG~\cite{white-dmrg-1992} algorithms for finding the ground state of a quantum spin chain.
 Similar formats were applied in statistics and can be traced back to works of Hitchcock in 1930s~\cite{hitchcock-rank-1927,hitchcock-sum-1927} and Tucker in 1960s~\cite{Tucker}.
 Tensor formats in three dimensions were rigorously studied in the numerical linear algebra community as a generalisation of the low--rank decomposition of matrices~\cite{ost-tucker-2008,fkst-chem-2008,gst-piers-2009,sav-bss-2010}.
 Generalisations to higher dimensions eventually led to re-discovery of the MPS/DMRG framework under the name of the tensor train (TT) format~\cite{ot-tt-2009,osel-tt-2011}.
 In the last decade the TT format and the more general Hierarchical Tucker format~\cite{gras-hsvd-2010} have been successfully applied to a variety of problems, such as
  finding several low--lying states of a quantum system~\cite{dkos-eigb-2014},
  performing superfast Fourier transform~\cite{dks-ttfft-2012},
  solving high--dimensional linear systems~\cite{ds-amen-2014} and
  first--order differential equations~\cite{d-tamen-2018}.
 First applications of tensor product methods to optimal control problems are considered in~\cite{dpss-frac-2015,bdos-sb-2016}.

 The recently proposed tAMEn algorithm~\cite{d-tamen-2018} allows to compute evolution of a quantum system under a given (time--dependent) Hamiltonian with high accuracy.
 Although the trajectory lies in extremely high--dimensional state space, tAMEn keeps the computation costs feasible by performing all calculations in the TT format and adapting the TT ranks to meet the required accuracy.
 In this paper we employ the tAMEn algorithm as a building block to develop a version of a classical GRAPE algorithm~\cite{grape-2005} and use it to control a spin chain with $d=41$ spins.
 This shows that tensor product algorithms can be used to design control sequences for quantum computers with $50$ to $100$ qubits, which can be expected within a decade~\cite{blackett-2016}.

 \newpage
\section{Optimal quantum control}                           
\subsection{Dynamic optimisation problem}
The dynamics of a quantum system with Hamiltonian $\H$ is described by Liouville--von Neumann equation in superoperator form
\begin{equation}\nonumber
 \dot\rho   = -\i   \hH (t) \rho  + \hR\rho,
\end{equation}
where $\rho$ is the quantum density matrix stretched into a vector form (referred to later as quantum state), and
\begin{equation}\nonumber
 \hH(t)   = \Id \otimes \H(t)  - \H(t) \otimes \Id.
\end{equation}
 Traditionally $\H(t)$ is split into a constant drift and a time--dependent control term:
\begin{equation}\nonumber
 \H(t)  = \H_0 + \sum_{k=1}^K c_k(t) \H_k.
\end{equation}
Control operators $\H_k$ are usually defined by model or instrument used; we manipulate the system by changing control pulse functions $c_k(t).$
Now choosing a suitable basis set $\{\phi_n(t)\}_{n=1}^N$ and dropping the relaxation part $\hR$ for simplicity, we obtain an ODE for system dynamics:
\begin{equation}\label{eq:ode}
 \dot\rho = -\i \left( \hH_0 + \sum_{k=1}^K \sum_{n=1}^{N} c_{n,k} \hH_k \phi_n(t) \right) \rho, 
 \qquad
 \rho(0)=\rho_0.
\end{equation}
Our goal is to choose control parameters $c=[c_{n,k}]$ to maximise \emph{fidelity} (or \emph{overlap})
\begin{equation}\label{eq:fidelity}
 F(c) = \Re \aprod{\rho_T}{\rho(T)}.
\end{equation}

\subsection{First--order optimisation framework}
Classical optimisation methods require a gradient of $F$ w.r.t. $c_{n,k}$ which reads
\begin{equation}\label{eq:grad}
\frac{\d F}{\d c_{n,k}} = \Re\aprod{\rho_T}{\frac{\d \rho(T)}{\d c_{n,k}}}.
\end{equation}
 To find the gradient, we differentiate \eqref{eq:ode} w.r.t. control parameters and obtain
\begin{equation}\nonumber
 \frac{\d}{\d c_{n,k}} \frac{\d \rho}{\d t}
   = -\i\left( \phi_n(t) \hH_k \rho + \hH(t) \frac{\d\rho}{\d c_{n,k}}\right).
\end{equation}
Changing the order of derivatives in the left hand side, we obtain a system of coupled ODEs for the density matrix and its gradient:
\begin{equation}\label{eq:ode2}
 \frac{\d}{\d t} 
   \begin{bmatrix}
      \frac{\d\rho}{\d c_{n,k}} \\ 
      \rho
   \end{bmatrix} 
   = 
   -\i
   \begin{bmatrix}
      \hH(t) & \phi_n(t)\hH_k \\ 
      0      & \hH(t)
   \end{bmatrix} 
   \: 
   \begin{bmatrix}
   \frac{\d\rho}{\d c_{n,k}} 
   \\ 
   \rho
 \end{bmatrix},
 \qquad
   \begin{bmatrix}
      \frac{\d\rho}{\d c_{n,k}} \\ 
      \rho
   \end{bmatrix}_{t=0}
   = 
   \begin{bmatrix}
      0           \\ 
      \rho_0
   \end{bmatrix}.
\end{equation}
To justify the initial condition for the gradient, write a first--order approximation for $\rho(\delta t)$ with an infinitesimal time $\delta t$:
$$
 \rho(\delta t) = \rho_0 - \i \delta t \left(\hH_0 + \sum_{k=1}^K \sum_{n=1}^{N} c_{n,k} \hH_k \phi_n(t)\right) \rho_0 + \O(\delta t^2)
$$
and differentiate it w.r.t. $c_{n,k}$ to see $\frac{\d \rho}{\d c_{n,k}}(\delta t) = \O(\delta t)\to0.$

The system \eqref{eq:ode2} can be integrated on $[0,T]$ using a suitable method (e.g. Runge--Kutta scheme~\cite{runge}) for all $k,n$ to produce both the fidelity~\eqref{eq:fidelity} and its gradient~\eqref{eq:grad}.

\subsection{GRAPE algorithm}
A classical gradient ascend pulse engineering (GRAPE) method~\cite{grape-2005} is a simple version of the approach described above.
Control pulses $c_k(t)$ are assumed piecewise--constant on intervals $(t_{n-1},t_n]$ between the nodes of the grid $\{t_n\}_{n=0}^N$ with $t_0=0$ and $t_N=T.$
We will use a uniform grid $t_k=\tau k$ with step--size $\tau=T/N.$
The basis functions $\phi_n(t)$ can be taken as Heaviside functions
\begin{equation}\nonumber
 \theta_n(t) = \begin{cases}
   1, & \text{if~ $t_{n-1} < t \leq t_n$},  \\ 
   0, & \text{otherwise},
 \end{cases}
\end{equation}
thus making $c_{n,k}$ the amplitude of control pulse $c_k(t)$ on $(t_{n-1},t_n].$
The time--locality of basis functions $\theta_n(t)$ allows us to compute~\eqref{eq:grad} more efficiently.
Note that $\theta_n(t)=0$ for $t\le t_{n-1}$ hence ODE~\eqref{eq:ode2} for the gradient component simplifies to the homogeneous ODE
\begin{equation} \label{eq:grad0}
 \frac{\d}{\d t}\frac{\d\rho}{\d c_{n,k}} = -\i \hH \frac{\d\rho}{\d c_{n,k}} 
 \quad\text{for}\quad t\in(0,t_{n-1}],
 \quad\text{s.t.}\quad \frac{\d\rho}{\d c_{n,k}}(0)=0,
\end{equation}
 which of course gives
 $
 \frac{\d\rho}{\d c_{n,k}} (t_{n-1}) = 0.
 $
 Therefore, to compute~\eqref{eq:grad} we can propagate the state $\rho(t)$ on $(0,t_{n-1}]$ and start solving the coupled system \eqref{eq:ode2} only from $t_{n-1}.$
 When we reach $t_n$ and compute
 $
 \frac{\d\rho}{\d c_{n,k}} (t_{n}),
 $
 we can note that $\theta_n(t)=0$ for $t>t_n$ hence \eqref{eq:grad0} holds on $t\in(t_n,T]$ as well.
 Since $\hH(t)$ is piecewise--constant,
\begin{equation}\nonumber
\frac{\d\rho}{\d c_{n,k}}(T) 
   = \exp \left[-\i\int_{t_n}^T \hH(s) \ds \right]\frac{\d\rho}{\d c_{n,k}}(t_n)  
   = \prod_{m={n+1}}^{N} \underbrace{\exp \left[-\i\tau \hH(t_{m})\right]}_{\Q_m} \frac{\d\rho}{\d c_{n,k}}(t_n),
\end{equation}
and plugging this in~\eqref{eq:grad} gives
\begin{equation}\label{eq:grad2}
 \frac{\d F}{\d c_{n,k}} 
  = \Re\aprod{\rho_T}{\prod_{m={n+1}}^{N} \Q_m \frac{\d\rho}{\d c_{n,k}}(t_n)} 
  = \Re\aprod{\prod_{m={n+1}}^{N} \Q^{\dagger}_m\rho_T}{\frac{\d\rho}{\d c_{n,k}}(t_n)}.
\end{equation}
The state 
$$
\lam(t_n)=\prod_{m={n+1}}^{N} \Q^{\dagger}_m\rho_T
$$
is the solution of ODE with the boundary condition at the right end, colloquially referred to as `back--propagation':
\begin{equation}\label{eq:ode-lambda}
 \frac{\d\lam}{\d t} = -\i \hH(t) \lam
 \quad\text{for}\quad t\in[t_{n},T),
 \quad\text{s.t.}\quad \lam(T) = \rho_T.
\end{equation}

Hence the coupled system \eqref{eq:ode2} needs to be solved only on $[t_{n-1}, t_n]$ to produce the $n$-th components of the gradient.
By pre-computing all $\lam(t_n)$ for $n=1,\ldots,N,$ we can implement the GRAPE iteration with complexity scaling linearly in $N.$
\begin{enumerate}
 \item For $n=N,\ldots,1,$ propagate \eqref{eq:ode-lambda} backward: \hfill 
  $
  \lam(t_{n}) \:\longrightarrow\: \lam(t_{n-1})
  $
 \item For $n=1,\ldots,N,$ propagate \eqref{eq:ode2} forward: \hfill
  $
  \begin{bmatrix}0 \\ \rho(t_{n-1})\end{bmatrix}
   \:\longrightarrow\:
  \begin{bmatrix}\frac{\d \rho}{\d c_{n,k}}(t_n) \\ \rho(t_n)\end{bmatrix}
  $
 \item Compute \eqref{eq:grad2} and update the control amplitudes: \hfill
  $
  c_{n,k} \:\longrightarrow\: c_{n,k} - \epsilon \frac{\d F}{\d c_{n,k}}.
  $
\end{enumerate}
Here $\epsilon$ is the step size for the gradient ascend.
It should be carefully selected, since a too small value will make the algorithm converging very slowly while a very big one will make the algorithm divergent.
In our implementation the step size is chosen adaptively:
 it is increased each time the step is successful (i.e. results in better fidelity) and 
 decreased each time we have an overstep (i.e. fidelity decreases and the step has to be rejected).
Note that line--search strategies for this optimisation problem are rather pointless, because the gradient~\eqref{eq:grad2} and the fidelity function~\eqref{eq:fidelity} are computed simultaneously from~\eqref{eq:ode2}, and have therefore the same complexity.

In the classic GRAPE algorithm, the propagators $\Q_m$ in~\eqref{eq:grad2} are  calculated as matrix exponentials~\cite{moler-expm-2003} and applied to the states;
in more advanced versions of the algorithm the action of matrix exponentials on the states is computed on--fly maintaining the sparse format for $\lam$'s and $\rho$'s~\cite{sidje-expokit-1998}.

Instead of computing propagators, we use the tAMEn algorithm to solve the first--order ODEs~\eqref{eq:ode2} and~\eqref{eq:ode-lambda} with Hamiltonians and states represented in the TT format.
The next section focuses on mathematical details of this approach.

\section{Tensor train format and the tAMEn algorithm}                
Introducing individual state indices for each site, we can treat $\rho$ as a multi--index (high--dimensional) array which is referred to as \emph{tensor} in numerical linear algebra:
$$
\rho=[\rho(i_1,\ldots,i_d;\: j_1,\ldots,j_d)].
$$
We should admit that our use of the term \emph{tensor} does not imply proper differentiation of upper and lower indices or \emph{Einstein summation convention}~\cite{Einstein-relativitat-1916}.
We use the Tensor Train (TT) decomposition \cite{osel-tt-2011} to separate \emph{pairs} of $i_k,j_k$, belonging to different sites, by a product approximation of the form
\begin{equation}\label{eq:tt}
 \rho(i_1,\ldots,i_d;~j_1,\ldots,j_d) \approx \sum_{\alpha_0,\ldots,\alpha_d=1}^{r_0,\ldots,r_d} \rho^{(1)}_{\alpha_0,\alpha_1}(i_1, j_1) \cdots \rho^{(d)}_{\alpha_{d-1},\alpha_d}(i_d, j_d).
\end{equation}
Here, $\rho^{(k)}$, $k=1,\ldots,d,$ are called \emph{TT blocks}, and the summation ranges $r_1,\ldots,r_d$ are called \emph{TT ranks}.
The grouping of $i_k$ and $j_k$ is such that unentangled (e.g. pure) states are represented by \eqref{eq:tt} exactly with elementary TT ranks $r_0=\cdots=r_d=1$.
The TT representation is equivalent to matrix product states (MPS)~\cite{fannes-mps-1992,klumper-mps-1993} with open boundary conditions $r_0=r_d=1$.
Entangled states require $r_k>1$ for at least one intermediate $k=1,\ldots,d-1.$
However, we aim at representing states with \emph{weak} correlations,
so that $r_1,\ldots,r_{d-1}$ can be kept bounded by a moderate constant $r_k \leq r \ll 2^d.$
This yields an $\O(dr^2)$ storage cost of the TT format, \emph{linear} in the system size.
For example, ground states of non--critical one--dimensional spin chains admit such bounded TT ranks due to the area law~\cite{Eisert-area-laws-2010}.

We use a linear piecewise Chebysh\"ev scheme~\cite{d-tamen-2018} to discretise ODEs~\eqref{eq:ode}, \eqref{eq:ode2} and \eqref{eq:ode-lambda} in time.
Recall that we solve ODEs on time intervals $t\in(t_{n-1}, t_n]$, where control pulses $c_k(t)$ and hence $\hH(t)$ are constant and only $\rho(t)$ depends on time.
Without loss of generality we can assume that $t\in(0,\tau]$ and drop the time index $n$ for the rest of this section.
We choose a basis of Lagrange polynomials $\{L_m(t)\}_{m=1}^M$ centered at Chebysh\"ev nodes $\{\tau_m\}_{m=1}^M$ and represent the state as
$$
\rho(t) \approx \sum_{m=1}^{M} \rho(\tau_m) L_m(t), 
\qquad t\in(0,\tau],
\qquad 0<\tau_1<\tau_2<\cdots<\tau_M = \tau.
$$ 
Since the length of the interval $\tau$ is usually quite small to provide better resolution of control pulses, the basis size $M$ can be very moderate.
We collect the nodal values $\rho(\tau_m)$ into a vector
$
\overline\rho = [\rho(\tau_m)]_{m=1}^{M}
$ 
of length $4^d M,$ which we aim to find.
The time derivative operator is discretised as a \emph{differentiation matrix}
$
S = [L'_{\ell}(\tau_m)]_{m,\ell=1}^M
$
in accordance to the spectral approximation theory~\cite{trefethen-spectral-2000}.
Applying this discretisation scheme to ODE~\eqref{eq:ode} we obtain the following linear system
\begin{equation}\label{eq:timescheme}
 (\underbrace{S \otimes \Id + \Id \otimes \i \hH}_{A} ) \overline\rho
 = 
 \underbrace{(S \: \e_M) \otimes \rho(0) \strut}_{f}.
\end{equation}
Here 
  $\overline\rho$ is the unknown vector of $\rho(\tau_m)$'s, 
  $\e_M$ is the vector of size $M$ full of $1$'s, and
  $\Id$ denotes identity matrices of appropriate size.
Note that the Kronecker products are implicitly realised by the TT representation~\eqref{eq:ttt}.
Therefore, we never actually compute them in \eqref{eq:timescheme}, but just collect the corresponding factors into TT blocks.
The same applies to the states.
In particular, we consider unentangled states $\rho_0$ and $\rho_T,$ and construct their TT blocks explicitly.
Simple mixed states can be made by concatenation of TT blocks instead of a direct summation of a large number of density matrices~\cite{sdwk-nmr-2014}.
The TT ranks can be truncated down to optimal values for a desired accuracy using the singular value decomposition (SVD)~\cite{osel-tt-2011}.

Algebraic equations \eqref{eq:timescheme} can be solved in the TT approximation by an alternating iteration,
such as the Alternating Minimal Energy (AMEn) algorithm \cite{ds-amen-2014},
which is an enhanced version of Alternating Least Squares (ALS) \cite{holtz-ALS-DMRG-2012}.
We expand the TT representation \eqref{eq:tt} by an extra block which carries the time index $m$ and write:
\begin{equation}\label{eq:ttt}
 \begin{split}
  \overline\rho(i_1,\ldots,i_d; j_1,\ldots,j_d,m) & = \sum_{\alpha} \rho^{(1)}_{\alpha_0,\alpha_1}(i_1,j_1) \cdots \rho^{(d)}_{\alpha_{d-1},\alpha_d}(i_d,j_d) \rho^{(d+1)}_{\alpha_d,\alpha_{d+1}}(m),
  \\
  \overline\rho  & = \sum_{\alpha} \rho^{(1)}_{\alpha_0,\alpha_1}\otimes \cdots \otimes \rho^{(d)}_{\alpha_{d-1},\alpha_d} \otimes \rho^{(d+1)}_{\alpha_d,\alpha_{d+1}}.
 \end{split}
\end{equation}
 Now the last block $\rho^{(d+1)} = [\rho^{(d+1)}_{\alpha_d,\alpha_{d+1}}(m)]$ carries degrees of freedom in time,
  the rightmost TT rank is $r_{d+1}=1,$ and $r_d$ can now be larger than $1.$

The basic ALS algorithm would solve~\eqref{eq:timescheme} by sweeping over TT blocks, i.e. representing $\overline\rho$ by \eqref{eq:ttt} and reducing the equation to only one block $\rho^{(k)}$ in each step.
This can be written using a \emph{frame} matrix
$$
P_{\neq k} 
 = \left(\sum_{\alpha_0\ldots\alpha_{k-2}} \rho^{(1)}_{\alpha_0,\alpha_1}\otimes \cdots \otimes \rho^{(k-1)}_{\alpha_{k-2}}\right)
        \otimes I \otimes 
   \left(\sum_{\alpha_{k+1}\ldots\alpha_{d+1}} \rho^{(k+1)}_{\alpha_{k+1}}\otimes \cdots \otimes \rho^{(d+1)}_{\alpha_{d},\alpha_{d+1}}\right),
$$
which is of size $4^dM \times 4 r_{k-1} r_k$ for $k\leq d,$ and of size $4^dM \times r_{k-1} M r_k$ for $k=d+1.$
The frame matrix realises a linear map from the elements of a single TT block into the elements of the whole vector given \eqref{eq:ttt}, that is 
$
\overline\rho = P_{\neq k} \rho^{(k)},
$
assuming that elements of $\rho^{(k)}$ are stretched in a vector.
The ALS method solves the reduced Galerkin systems
$
(P_{\neq k}^{\dagger} A P_{\neq k})\rho^{(k)} = P_{\neq k}^{\dagger} f
$ 
subsequently for $k=1,\cdots,d+1.$
This system can be assembled efficiently~\cite{holtz-ALS-DMRG-2012} due to the TT formats of $P_{\neq k}$, $A$ and $f$ and solved using standard methods.

The AMEn method \cite{ds-amen-2014} improves the convergence of ALS by expanding $\rho^{(k)}$ (and hence all $P_{\neq q}$ for $q>k$) with a TT approximation of the \emph{residual} $f-A\overline\rho,$ where $\overline\rho$ is formed by~\eqref{eq:ttt} with an updated block $\rho^{(k)}$ plugged in.
This also allows to adapt the TT ranks to ensure a desired accuracy of the TT approximation.
The tAMEn (time--dependent AMEn) algorithm~\cite{d-tamen-2018} utilises the special meaning of the last TT block, carrying the time variable $m,$ in order to preserve conservation laws (e.g. the Frobenius norm), and to estimate the time discretisation error.

For the interval lengths $\tau$ chosen in our numerical experiments,
we have found that $M=8$ Chebysh\"ev nodes in each interval are sufficient to resolve the time derivative with the relative accuracy of $10^{-5}$ or better.
In each step the tAMEn algorithm produces a set of TT blocks $\rho^{(1)},\ldots,\rho^{(d+1)},$ comprising the discrete--time solution $\overline\rho$ on the interval $(0,\tau].$
After $\overline\rho$ is obtained, we compute $\rho(\tau)$ in the TT format by taking the slice of the last TT block $\rho^{(d+1)}(M)$ and hence removing the time index $m$ from consideration.
Since $t=\tau$ represents the end point of the interval $(t_{n-1},t_n]$ the obtained state represents $\rho(t_n)$ and can be used now as an initial state in the next GRAPE step over the time interval $(t_n,t_{n+1}].$

The tAMEn algorithm is agnostic to a particular form of $\hH$ (the only assumption is that it must admit a TT decomposition), and can be applied to any tensor--structured differential equation $\dot x = Ax$
such as the auxiliary matrix formalism \eqref{eq:ode2} and the back--propagation \eqref{eq:ode-lambda}.

\section{Numerical experiments}                             
 We consider a Heisenberg chain of spin--$\frac12$ particles, for which a Hamiltonian is a simple sum of nearest neighbour interactions
 \begin{equation}\label{eq:H}
 \H = \sum_{k=1}^{d-1} 
         J_x \s_x^{(k)} \s_x^{(k+1)} 
       + J_y \s_y^{(k)} \s_y^{(k+1)} 
       + J_z \s_z^{(k)} \s_z^{(k+1)},
 \end{equation}
 where the Pauli matrices $\s_{\{x,y,z\}}^{(k)}$~\cite{slavnov-bethe-2018} act only on spin in position $k$ of the chain,
 $$
 \s_{\{x,y,z\}}^{(k)} = \Id \otimes \cdots \otimes \Id \otimes \underbrace{\s_{\{x,y,z\}}}_{\text{$k$--th place}} \otimes \Id \otimes \ldots \otimes \Id.
 $$
 When $J_x=J_y=J_z,$ this system is called the XXX Heisenberg model, when $J_x=J_y\neq J_z$ it is called the XXZ model, and for $J_x\neq J_y \neq J_z$ it is called the XYZ model.
 It is commonly mentioned that linear Heisenberg chains can be diagonalised exactly using the Bethe ansatz~\cite{bethe-1931}, 
   however the eigenstates can only be written via roots of a system of algebraic equations of degree $d$~\cite{slavnov-bethe-2018},
   that are not computable in closed form neither using analytic methods nor numerically with reasonable accuracy for $d\gtrsim 20.$
 For simpler XX models ($J_x=J_y,$ $J_z=0$) eigenvectors are available in closed form and can be used for dimensionality reduction, allowing chains with $d \simeq 200$ spins to be controlled~\cite{spinwire-simple-2005}.
 We could not find examples of optimal control pulse computed numerically for XXX, XXZ and XYZ Heisenberg models with the number of spins $d\gtrsim 20.$
 We decided therefore to test our tensor product approach for XXX and XXZ Heisenberg spin-$\frac12$ chains with $d=11,$ $d=21$ and $d=41.$

 The initial and target states are taken as
 \begin{equation}\nonumber
  \psi_0 = \ket{\uparrow \downarrow \downarrow \cdots \downarrow \downarrow}, \qquad
  \psi_T = \ket{\downarrow \downarrow \cdots \downarrow \downarrow \uparrow},
 \end{equation}
 so the task is to move the $\ket{\uparrow}$ state from the first to the last position in the chain.
 The control operator  $H_c(t) = c(t) \s_z^{(1)}$ is the magnetic field acting on the first spin only.
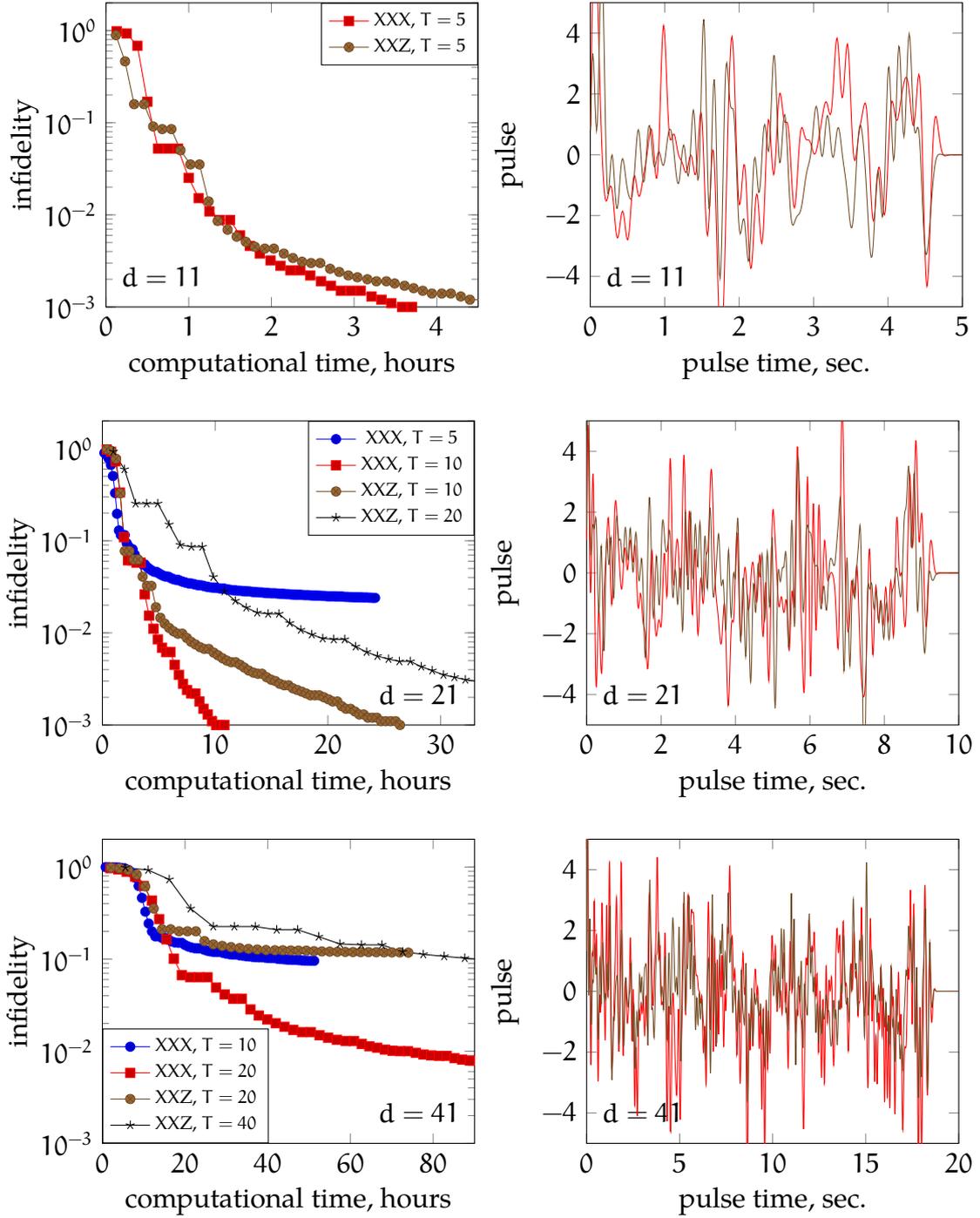
\begin{figure}[p]
 \begin{center}
  \begin{tikzpicture}
   \begin{axis}[width=.48\textwidth,
     xmode=normal,  xmin=0,    xmax=4.5,
     ymode=log,     ymin=1e-3, ymax=2,
     xlabel={computational time, hours},
     ylabel={infidelity},
     legend style={at={(1,1)},anchor=north east,draw=black,thin,fill=white},
     ]
     \node at(axis description cs:0.15,0.1){$d=11$};
     \pgfplotsset{cycle list shift=1};
     \addplot+[] table[header=true,x=hours,y=fidelity]{./dat/xxx_d11_t05.dat}; \addlegendentry{\scriptsize XXX, $T=5$};
     \addplot+[] table[header=true,x=hours,y=fidelity]{./dat/xxz_d11_t05.dat}; \addlegendentry{\scriptsize XXZ, $T=5$};
   \end{axis}
  \end{tikzpicture}
  \begin{tikzpicture}
   \begin{axis}[width=.48\textwidth,
     xmode=normal, xmin=0,  xmax=5,
     ymode=normal, ymin=-5, ymax=5,
     xlabel={pulse time, sec.},
     ylabel={pulse},
     legend style={at={(1,1)},anchor=north east,draw=black,thin,fill=white},
     grid=none
     ]
     \node at(axis description cs:0.15,0.1){$d=11$};
     \pgfplotsset{cycle list shift=1};
     \addplot+[no marks] table[header=false,x index=0,y index=1]{./dat/xxx_d11_t05.out};
     \addplot+[no marks] table[header=false,x index=0,y index=1]{./dat/xxz_d11_t05.out};
   \end{axis}
  \end{tikzpicture}
 \end{center}
 \begin{center}
  \begin{tikzpicture}
   \begin{axis}[width=.48\textwidth,
     xmode=normal, xmin=0,    xmax=33,
     ymode=log,    ymin=1e-3, ymax=2,
     xlabel={computational time, hours},
     ylabel={infidelity},
     legend style={at={(1,1)},anchor=north east,draw=black,thin,fill=white},
     ]
     \node at(axis description cs:0.85,0.1){$d=21$};
     \addplot+[] table[header=true,x=hours,y=fidelity]{./dat/xxx_d21_t05.dat}; \addlegendentry{\scriptsize XXX, $T=5$};
     \addplot+[] table[header=true,x=hours,y=fidelity]{./dat/xxx_d21_t10.dat}; \addlegendentry{\scriptsize XXX, $T=10$};
     \addplot+[] table[header=true,x=hours,y=fidelity]{./dat/xxz_d21_t10.dat}; \addlegendentry{\scriptsize XXZ, $T=10$};
     \addplot+[] table[header=true,x=hours,y=fidelity]{./dat/xxz_d21_t20.dat}; \addlegendentry{\scriptsize XXZ, $T=20$};
   \end{axis}
  \end{tikzpicture}
  \begin{tikzpicture}
   \begin{axis}[width=.48\textwidth,
     xmode=normal, xmin=0,  xmax=10,
     ymode=normal, ymin=-5, ymax=5,
     xlabel={pulse time, sec.},
     ylabel={pulse},
     legend style={at={(1,1)},anchor=north east,draw=black,thin,fill=white},
     grid=none
     ]
     \node at(axis description cs:0.15,0.1){$d=21$};
     \pgfplotsset{cycle list shift=1};
     \addplot+[no marks] table[header=false,x index=0,y index=1]{./dat/xxx_d21_t10.out};
     \addplot+[no marks] table[header=false,x index=0,y index=1]{./dat/xxz_d21_t10.out};
   \end{axis}
  \end{tikzpicture}
 \end{center}
 \begin{center}
  \begin{tikzpicture}
   \begin{axis}[width=.48\textwidth,
     xmode=normal, xmin=0,    xmax=90,
     ymode=log,    ymin=1e-3, ymax=2,
     xlabel={computational time, hours},
     ylabel={infidelity},
     legend style={at={(0,0)},anchor=south west,draw=black,thin,fill=white},
     ]
     \node at(axis description cs:0.85,0.1){$d=41$};
     \addplot+[] table[header=true,x=hours,y=fidelity]{./dat/xxx_d41_t10.dat}; \addlegendentry{\scriptsize XXX, $T=10$};
     \addplot+[] table[header=true,x=hours,y=fidelity]{./dat/xxx_d41_t20.dat}; \addlegendentry{\scriptsize XXX, $T=20$};
     \addplot+[] table[header=true,x=hours,y=fidelity]{./dat/xxz_d41_t20.dat}; \addlegendentry{\scriptsize XXZ, $T=20$};
     \addplot+[] table[header=true,x=hours,y=fidelity]{./dat/xxz_d41_t40.dat}; \addlegendentry{\scriptsize XXZ, $T=40$};
   \end{axis}
  \end{tikzpicture}
  \begin{tikzpicture}
   \begin{axis}[width=.48\textwidth,
     xmode=normal, xmin=0,  xmax=20,
     ymode=normal, ymin=-5, ymax=5,
     xlabel={pulse time, sec.},
     ylabel={pulse},
     legend style={at={(1,1)},anchor=north east,draw=black,thin,fill=white},
     grid=none
     ]
     \node at(axis description cs:0.15,0.1){$d=41$};
     \pgfplotsset{cycle list shift=1};
     \addplot+[no marks] table[header=false,x index=0,y index=1]{./dat/xxx_d41_t20.out};
     \addplot+[no marks] table[header=false,x index=0,y index=1]{./dat/xxz_d41_t20.out};
   \end{axis}
  \end{tikzpicture}
 \end{center}
 \caption{Pulse sequence optimisation for the Heisenberg XXX and XXZ spin chains given by~\eqref{eq:H} with $J_x=J_y=J_z=2\pi$ and  $J_x=J_y=2\pi, J_z=2.2\pi,$ respectively.
 Left: convergence of the GRAPE algorithm with TT compression of all states and tAMEn algorithm for time evolution.
 Right: the optimised pulse sequences for XXX and XXZ chains. 
 Top, middle, bottom row: chains with $d=11,$ $d=21$ and $d=41$ spins, respectively.
 }
 \label{fig}
\end{figure}

 In theory, Heisenberg chains are fully controllable~\cite{schirmer-control-2016}, which means the fidelity~\eqref{eq:fidelity} can  be made infinitely close to $1$ as $T\to\infty.$
 In practice, however, the available time $T$ for the pulse is limited, and the final state $\rho(T)$ will not be fully focused, leaving some \emph{infidelity}
 $
  1 - | \aprod{\rho_T}{\rho(T)}| > 0.
 $
 We optimise the pulse sequence to reduce the infidelity as much as possible.

 The results are shown on~Fig.~\ref{fig}.
 Computing a short high--fidelity pulse for large spin system appears to be a challenging task.
 Using 4 cores on a Xeon E5-2650 CPU and 10 GB of memory, we reduced infidelity to $10^{-3}$ for $d=21$ and to $10^{-2}$ for $d=41$ within a few days of computation.
 This seems fast compared to high--fidelity calculations reported in~\cite{spiteri-highfidelity-control-2018}, where infidelity $10^{-4}$ was reached after weeks/month of optimisation for a system with $d=4$ qubits with $4$ levels each.

 The efficiency of the method depends primarily on reasonable choice of the pulse length $T,$ which is usually of the order of the inverse of the natural frequency of the system~\cite{khaneja-control-2001}. 
 The time $T$ typically should increase with the number of spins in the chain due to the finite rate of the information exchange among them.
 We see in our experiments that larger intervals $T$ are required to reduce infidelity for larger chains; also, XXZ chains demand larger $T$ than XXX chains.
 We keep resolution $\tau=10^{-2}$ for the control pulse constant for all experiments, which is the main reason why computational time per iteration grows with $T.$
 The TT ranks grow only mildly with $d$ and remain moderate $r_k\lesssim 30$ for $d=41.$
 This resulted in about $70\cdot10^3$ unknowns in the TT decomposition of each $\rho(t_k)$ for $d=21$ and $140\cdot10^3$ unknowns for $d=41$.
 For comparison, a single diagonal of $\rho(t_k)$ has $2^{21}\approx 2\cdot10^6$ entries for $d=21$ and $2^{41}=2\cdot10^{12}$ entries for $d=41.$

\newpage
\section{Conclusions and future work}                                       
Tensor product algorithms open new possibilities to control long spin wires and multi--qubit gates in quantum computers.
We explored a proof-of-concept example, in which 
 we considered a simple XXX or XXZ Heisenberg chain, 
 optimised control pulse using a classical GRAPE method~\cite{grape-2005}, while 
 representing all states n the TT format~\cite{osel-tt-2011} and propagating them using the tAMEn algorithm~\cite{d-tamen-2018}.
This combination of relatively simple algorithms allowed us to reach fidelity of $99\%$ for a chain of $d=41$ spins using a single workstation for calculations.
We used only controllable approximation techniques and avoided any heuristic or random truncations of the state space.
Our algorithm is deterministic and flexible, i.e. can be applied to any linear or quasi--linear quantum system, for which we can expect the states to have moderate entanglement (as measured by Schmidt ranks).

We have not  yet reached the fidelity of $\sim99.99\%$ required for topological error correction~\cite{spiteri-highfidelity-control-2018}, mostly because of slow first--order convergence of GRAPE, but we hope to make it possible using a second--order optimisation algorithm, such as Newton--Raphson~\cite{kuprov-grape-2016} or a quasi--Newton e.g. BFGS~\cite{kuprov-grape-bfgs-2011,eitan-grape-bfgs-2011}.
We have not imposed any constraints on pulse shape, which resulted in a few outbreaks seen on Fig.~\ref{fig}, but implementation of the box constraints $|c(t)|<C$ is a straightforward step which we postpone for future work.

\section*{Acknowledgements}
 Authors appreciate financial support of EPSRC which makes this work possible --- grants EP/P033954/1 (D.Q. and D.S.) and EP/M019004/1 (S.D.)
 Numerical calculations for this paper were performed on a dedicated server, provided to D.S. by the University of Brighton (Rising Stars grant).

 \newpage

\end{document}